\documentclass[11pt,reqno]{amsart}
\usepackage{geometry}
\usepackage[T1]{fontenc}
\usepackage[utf8]{inputenc}
\usepackage{lmodern}
\usepackage{amsmath,amssymb,mathtools}
\usepackage{mathrsfs}
\usepackage{microtype}
\usepackage[colorlinks=true, urlcolor=blue,  linkcolor=blue, citecolor=blue]{hyperref}
\usepackage{enumitem}
\newtheorem{theorem}{Theorem}[section]
\newtheorem{proposition}[theorem]{Proposition}
\newtheorem{lemma}[theorem]{Lemma}
\newtheorem{corollary}[theorem]{Corollary}
\theoremstyle{definition}
\newtheorem{definition}[theorem]{Definition}
\newtheorem{example}[theorem]{Example}
\theoremstyle{remark}
\newtheorem{remark}[theorem]{Remark}

\DeclareMathOperator{\dist}{dist}
\DeclareMathOperator{\ran}{ran}
\DeclareMathOperator{\Der}{Der}
\DeclareMathOperator{\CB}{CB}
\DeclareMathOperator{\Bil}{Bil}
\DeclareMathOperator{\diag}{diag}
\newcommand{\A}{\mathcal A}
\newcommand{\J}{\mathcal J}
\newcommand{\B}{\mathcal B}
\newcommand{\C}{\mathbb C}

\newcommand{\norm}[1]{\left\lVert #1\right\rVert}
\newcommand{\cb}{\mathrm{cb}}
\newcommand{\mat}{\mathrm{mat}}

\numberwithin{equation}{section}
\allowdisplaybreaks[1]
\setlist[enumerate]{label=\textup{(\roman*)},leftmargin=2.1em}
\title[Matricial stability of generalized Jordan derivations]{Matricial stability and central corrections for generalized Jordan derivations on amenable $C^*$-algebras}

\author[M. Rahimi-Alangi]{Mohammad Rahimi-Alangi}
\address[M. Rahimi-Alangi]{Payame Noor University, Theran, Iran}
\email{mohammadalangi@pnu.ac.ir}

\author[S. Hashemi Sababe]{Saeed Hashemi Sababe$^*$}
\address[S. Hashemi Sababe]{Independent Senior Researcher, Edmonton, Alberta, Canada}
\email{hashemi\_1365@yahoo.com, sa.hashemi@ualberta.ca}
\thanks{$^*$Corresponding author}

\subjclass[2020]{Primary 46L05; Secondary 47B47}
\keywords{Generalized Jordan derivation, quantitative rigidity, operator bimodule, approximate diagonal, completely bounded map}

\begin{document}
\begin{abstract}
We study norm stability of third-type generalized Jordan identities, allowing their auxiliary maps to be arbitrary. A canonical normalization eliminates these witnesses with a factor of four. For a unital $C^*$-algebra and a unital contractive Banach bimodule, the zero-defect maps are precisely sums of a derivation and multiplication by a central module element. We characterize scalar stability by the closed range of the normalized Jordan-defect operator. For operator bimodules, two explicit $2\times2$ tests recover an ordered Leibniz defect and a commutator defect. If the algebra has an approximate diagonal bounded by $M$ and the module is weak$^*$ closed, these tests yield approximation by an inner derivation plus a central multiplier with error at most $8M$ times the normalized matrix defect. The corresponding completely bounded estimate holds at all matrix levels. For full matrix algebras the constant is $8$, independently of their size, and weak$^*$ closedness is unnecessary. When the map vanishes at the identity, the constants improve to $2M$ and $2$, respectively.
\end{abstract}
\maketitle

\section{Introduction}

The study of derivations on $C^*$-algebras connects algebraic identities with continuity, module structure, and cohomology. Sakai proved automatic continuity for derivations from a $C^*$-algebra into itself \cite{Sakai}. Ringrose subsequently established automatic continuity for derivations into Banach bimodules \cite{Ringrose}. These results explain why an exact Leibniz identity has strong analytic consequences even when continuity is absent from its definition. Quantitative questions ask for additional information: if an identity has a small defect on a unit ball, how close is the map to an exact solution?

Jordan identities provide a particularly useful setting for this question. A Jordan derivation satisfies
\[
 D(a\circ b)=D(a)\circ b+a\circ D(b),\qquad
 a\circ b=\tfrac12(ab+ba).
\]
Johnson proved that a bounded Jordan derivation from a $C^*$-algebra into a Banach bimodule is an associative derivation \cite{JohnsonJordan}. His treatment also introduced symmetric amenability and used it to obtain Jordan rigidity for suitable Banach algebras. Peralta and Russo later obtained the corresponding automatic-continuity conclusion without initially assuming that the Jordan derivation is bounded \cite{PeraltaRusso}. Ozawa's work on Dixmier approximation and symmetric amenability supplies further structural context for the averaging methods associated with Jordan identities \cite{Ozawa}.

The immediate starting point of this article is the work of Hosseini, Peralta, and Su \cite{HPS}. They compare three notions of generalized derivation whose defining identities involve, respectively, an associated derivation, a fixed element of a bidual module, or two auxiliary maps. Their third-type identity has the form
\[
 F(ab)=G(a)b+aH(b),
\]
where $F$ is linear and the witnesses $G,H$ need not be linear or continuous. Their automatic-continuity and equivalence results settle the exact associative problem for $C^*$-algebra domains. They also study the Jordan version
\begin{equation}\label{eq:intro-third}
 F(a\circ b)=G(a)\circ b+a\circ H(b).
\end{equation}
In the arXiv version supplied with this study, Theorem~3.1 proves continuity and Proposition~3.2 gives an associative generalized-derivation identity when the bimodule is essential. Numbered references to their results below refer to arXiv:2403.18773v2.

The exact conclusion does not by itself give a norm estimate for approximate solutions of \eqref{eq:intro-third}. Two issues arise. First, the witnesses are unrestricted, so an estimate must be expressed in terms of the principal map rather than their norms. Second, the symmetrized product suppresses the order of multiplication. The value $F(1)$ creates an additional obstruction: an associative generalized derivation can have a noncentral correction, whereas the unital third-type Jordan identity forces that correction to be central. Proposition~\ref{prop:exact} makes this consequence explicit for every unital contractive Banach bimodule.

There is substantial earlier work on quantitative algebraic identities. Johnson's cohomological theory of amenable Banach algebras provides approximate diagonals and the basic mechanism for correcting small associative cocycle defects \cite{JohnsonCohomology,JohnsonDiagonals}. Haagerup proved that nuclear $C^*$-algebras are amenable \cite{Haagerup}. Samei and Soltani Farsani established explicit hyperreflexivity estimates for bounded cocycle spaces, under a Banach-space hypothesis on the appropriate Hochschild cohomology group \cite{SameiFarsani}. Their results show why a closed-range condition is relevant to estimates measuring distance from a cocycle space. Alaminos, Extremera, Godoy, and Villena obtained quantitative zero-product determination for $C^*$-algebras and other Banach algebras \cite{ZeroProduct}.

Stability of approximate generalized Jordan derivations has a separate functional-equation literature. Eshaghi Gordji and Ghobadipour study Hyers--Ulam--Rassias stability and superstability for an identity involving an associated Jordan derivation, together with control of additive and homogeneous defects \cite{NearlyJordan}. Park, Ghaffary Ghaleh, and Ghasemi treat approximate $n$-Jordan $*$-derivations \cite{Park}. In our setting the principal map is already bounded and complex linear, and the error is measured on two unit balls. The two witnesses in the third-type identity have no derivation or regularity requirements. The passage from these witnesses to a central correction, followed by estimates at matrix levels, is therefore the specific problem studied below.

Our argument combines three ingredients. Evaluating at the identity produces a canonical bounded linear witness and a defect comparable to the infimum over arbitrary witnesses. A quotient-space argument identifies the precise functional-analytic obstruction to scalar stability. Finally, two block-matrix substitutions recover both the unsymmetrized product identity and the commutator with the correction element. The last step allows the usual approximate-diagonal correction to be applied with explicit constants, without a quantitative version of Johnson's Jordan-to-associative theorem.

More precisely, let $\varepsilon_2(T)$ be the norm of the canonical Jordan defect of the entrywise amplification $T_2$, and let $M$ bound an approximate diagonal. For a weak$^*$-closed unital operator bimodule we prove that there are $\eta\in X$ and $z_0\in Z_{\A}(X)$ such that
\begin{equation}\label{eq:intro-bound}
 \norm{T-(a\mapsto a\eta-\eta a+az_0)}
 \leq 8M\varepsilon_2(T).
\end{equation}
The same estimate holds in completely bounded norm when the right side uses the supremum of the defects at all matrix levels. A full matrix algebra has an exact diagonal of projective norm one; hence its constants are independent of its dimension. The constants are explicit upper bounds, with no claim of optimality.

Section~2 fixes the module and matrix conventions. Section~3 establishes canonical normalization, exact centrality, and the scalar closed-range criterion. Section~4 proves the block estimates, the correction theorem, and its matrix-algebra consequences. Section~5 discusses the scalar and nonunital questions suggested by these results.

\section{Preliminaries}

All algebras, modules, and linear maps are complex. Throughout, $\A$ is a unital $C^*$-algebra with identity $1$. A \emph{unital contractive Banach $\A$-bimodule} is a Banach space $X$ with associative left and right actions such that
\[
 1x=x=x1,\qquad \norm{ax}\leq\norm a\norm x,
 \qquad \norm{xa}\leq\norm x\norm a.
\]
The left and right actions commute. We use
\[
 a\circ x=x\circ a=\tfrac12(ax+xa),\qquad
 Z_{\A}(X)=\{x\in X:ax=xa\text{ for all }a\in\A\}.
\]
For $z\in X$, put $L_z(a)=az$. Unitality implies $\norm{L_z}=\norm z$.

A derivation $d:\A\to X$ is a linear map satisfying
\[
 d(ab)=d(a)b+ad(b).
\]
Write $\Der(\A,X)$ for the bounded derivations and define
\begin{equation}\label{eq:exact-class}
 \J(\A,X)=\{d+L_z:d\in\Der(\A,X),\ z\in Z_{\A}(X)\}.
\end{equation}
Every derivation into a unital module satisfies $d(1)=0$. Consequently the representation in \eqref{eq:exact-class}, when it exists, is unique: $z=S(1)$ and $d=S-L_z$.

We write $\B(\A,X)$ for the Banach space of bounded linear maps and $\Bil_s(\A\times\A,X)$ for the bounded symmetric bilinear maps, with norm
\[
 \norm\Phi=\sup_{\norm a,\norm b\leq1}\norm{\Phi(a,b)}.
\]
Distances to subspaces are taken in operator norm unless a subscript $\cb$ is displayed.

For matrix estimates we use a concrete operator-bimodule convention. Let $H,K$ be Hilbert spaces and let
\[
 \sigma:\A\longrightarrow\B(K),\qquad
 \rho:\A\longrightarrow\B(H)
\]
be unital $*$-representations. A \emph{unital operator $\A$-bimodule} is a norm-closed subspace $X\subseteq\B(H,K)$ invariant under
\[
 ax=\sigma(a)x,\qquad xa=x\rho(a).
\]
We call it \emph{weak$^*$ closed} if it is closed for the canonical weak$^*$ topology of $\B(H,K)$. In that case its norm balls are weak$^*$ compact and its module actions are weak$^*$ continuous. This concrete convention suffices for all results below.

The space $M_n(X)$ carries its norm as a subspace of $\B(H^n,K^n)$. The amplified module actions are matrix multiplication with the amplified representations. In particular,
\begin{equation}\label{eq:matrix-actions}
 \norm{AU}\leq\norm A\norm U,\qquad
 \norm{UA}\leq\norm U\norm A
 \quad(A\in M_n(\A),\ U\in M_n(X)).
\end{equation}
Block-coordinate inclusions are isometric and block-coordinate compressions are contractive. For $T\in\B(\A,X)$, its amplification is $T_n([a_{ij}])=[T(a_{ij})]$. It is bounded for each fixed $n$. Indeed, if $A=[a_{ij}]$, $\xi=(\xi_j)\in H^n$, and $\zeta=(\zeta_i)\in K^n$ are unit vectors, then
\[
 \left|\sum_{i,j}\langle T(a_{ij})\xi_j,\zeta_i\rangle\right|
 \leq\norm T\norm A
 \left(\sum_j\norm{\xi_j}\right)
 \left(\sum_i\norm{\zeta_i}\right)
 \leq n\norm T\norm A.
\]
Thus $\norm{T_n}\leq n\norm T$. A map is completely bounded if
\[
 \norm T_{\cb}=\sup_{n\geq1}\norm{T_n}<\infty.
\]
The space of such maps is denoted by $\CB(\A,X)$, and
\[
 \J_{\cb}(\A,X)=\J(\A,X)\cap\CB(\A,X).
\]
For an operator bimodule, $\norm{L_z}_{\cb}=\norm z$. The inner derivation
\[
 d_\eta(a)=a\eta-\eta a
\]
is completely bounded and satisfies $\norm{d_\eta}_{\cb}\leq2\norm\eta$.

We recall the diagonal convention used in the correction argument. The Banach projective tensor product $\A\widehat\otimes_\pi\A$ has multiplication map $\pi(u\otimes v)=uv$ and module actions
\[
 a\cdot(u\otimes v)=au\otimes v,
 \qquad (u\otimes v)\cdot a=u\otimes va.
\]
For an algebraic tensor its projective norm is
\[
 \norm m_\pi=\inf\left\{\sum_j\norm{u_j}\norm{v_j}:
 m=\sum_j u_j\otimes v_j\right\}.
\]

\begin{definition}\label{def:diagonal}
A bounded approximate diagonal with bound $M$ is a net $(m_\lambda)$ in $\A\widehat\otimes_\pi\A$ such that
\[
 \sup_\lambda\norm{m_\lambda}_\pi\leq M,\qquad
 \norm{a\cdot m_\lambda-m_\lambda\cdot a}_\pi\longrightarrow0
 \quad(a\in\A),\qquad
 \norm{\pi(m_\lambda)-1}\longrightarrow0.
\]
An exact diagonal is a tensor $m$ satisfying $a\cdot m=m\cdot a$ and $\pi(m)=1$.
\end{definition}

The approximate-diagonal characterization of amenability is part of Johnson's theory \cite{JohnsonCohomology,JohnsonDiagonals}. In particular, Haagerup's theorem implies that every unital nuclear $C^*$-algebra has such a diagonal with some finite bound \cite{Haagerup}. Our estimates retain the chosen bound $M$ explicitly.

\section{Canonical normalization and the scalar stability obstruction}

For $T\in\B(\A,X)$, define its third-type Jordan defect by
\begin{equation}\label{eq:q-defect}
 q(T)=\inf_{G,H:\A\to X}
 \sup_{\norm a,\norm b\leq1}
 \norm{T(a\circ b)-G(a)\circ b-a\circ H(b)}.
\end{equation}
The infimum ranges over all maps, without linearity, continuity, or measurability assumptions. It is finite because one may take $G=H=0$.

Put $z=T(1)$ and define the canonical witness and defect by
\begin{align}
 K_T(a)&=T(a)-\tfrac12a\circ z
       =T(a)-\tfrac14(az+za),\label{eq:canonical}\\
 \Psi_T(a,b)&=T(a\circ b)-K_T(a)\circ b-a\circ K_T(b),\label{eq:psi}\\
 \varepsilon_1(T)&=\norm{\Psi_T}.
\end{align}
The bilinear map $\Psi_T$ is symmetric, and $K_T$ is bounded and linear.

\begin{lemma}\label{lem:witness}
For every bounded linear $T:\A\to X$,
\begin{equation}\label{eq:witness-bound}
 q(T)\leq\varepsilon_1(T)\leq4q(T).
\end{equation}
Moreover, $q$ is a seminorm, $|q(T)-q(S)|\leq\norm{T-S}$, and the linear operator
\begin{equation}\label{eq:scalar-R}
 \mathscr R_1:\B(\A,X)\longrightarrow\Bil_s(\A\times\A,X),
 \qquad \mathscr R_1(T)=\Psi_T
\end{equation}
has norm at most $4$.
\end{lemma}

\begin{proof}
Taking $G=H=K_T$ gives the first inequality in \eqref{eq:witness-bound}. To prove the second, fix $r>q(T)$ and choose witnesses $G,H$ whose residual on the two unit balls is at most $r$. Interchanging $a$ and $b$ in the residual and averaging shows that the single map $K=(G+H)/2$ satisfies
\[
 E(a,b):=T(a\circ b)-K(a)\circ b-a\circ K(b),
 \qquad \norm{E(a,b)}\leq r
\]
whenever $\norm a,\norm b\leq1$. This averaging uses the commutativity of the Jordan product and of the notation $a\circ x=x\circ a$; it requires no regularity of the witnesses.

Evaluation at the identity gives
\[
 E(a,1)=T(a)-K(a)-a\circ K(1),\qquad
 E(1,1)=z-2K(1).
\]
Consequently, for every $a$ in the unit ball,
\begin{align*}
 K_T(a)-K(a)
 &=T(a)-\tfrac12a\circ z-K(a)\\
 &=E(a,1)+a\circ\left(K(1)-\tfrac12z\right)\\
 &=E(a,1)-\tfrac12a\circ E(1,1).
\end{align*}
Contractivity of the module actions therefore gives
\[
 \norm{K_T(a)-K(a)}\leq\tfrac32r\qquad(\norm a\leq1).
\]
We use this estimate only on the unit ball: the map $K$ may be nonlinear. For $a,b$ in the unit ball,
\[
 \Psi_T(a,b)=E(a,b)
 -(K_T(a)-K(a))\circ b-a\circ(K_T(b)-K(b)),
\]
so $\norm{\Psi_T(a,b)}\leq r+3r=4r$. Letting $r$ decrease to $q(T)$ proves the second inequality.

Witnesses for two maps may be added, and witnesses may be multiplied by a scalar. Taking infima proves the triangle inequality and absolute homogeneity for $q$. For fixed witnesses for $S$, replacing $S$ by $T$ increases the residual by at most $\norm{T-S}$ because $\norm{a\circ b}\leq1$. Hence $q(T)\leq q(S)+\norm{T-S}$, and exchanging $S,T$ proves the Lipschitz estimate.

Finally, $T\mapsto T(1)$, $T\mapsto K_T$, and $T\mapsto\Psi_T$ are linear. Since $\norm z\leq\norm T$, we have $\norm{K_T}\leq\frac32\norm T$. Thus
\[
 \norm{\Psi_T}\leq\norm T+2\norm{K_T}\leq4\norm T,
\]
which proves the claimed operator bound.
\end{proof}

The next statement records a unital consequence of the exact theorem in the source manuscript. Its associative generalized-derivation input is explicitly attributed; the additional centrality argument is given here in the original module, rather than through an extension of the witnesses to a bidual.

\begin{proposition}\label{prop:exact}
Let $X$ be a unital contractive Banach $\A$-bimodule. For $T\in\B(\A,X)$ the following conditions are equivalent:
\begin{enumerate}
\item $T$ satisfies \eqref{eq:intro-third} for some arbitrary witnesses $G,H$;
\item $q(T)=0$;
\item $\Psi_T=0$;
\item $T=d+L_z$ for a bounded derivation $d$ and $z\in Z_{\A}(X)$.
\end{enumerate}
In condition \textup{(iv)}, necessarily $z=T(1)$ and $d=T-L_z$.
\end{proposition}

\begin{proof}
Condition (i) implies (ii), and Lemma~\ref{lem:witness} proves the equivalence of (ii) and (iii), even if the infimum defining $q$ is not attained. Condition (iii) implies (i) by taking $G=H=K_T$.

Suppose now that (iii) holds. By Proposition~3.2 of \cite{HPS}, there exists
$\xi\in X^{**}$ such that
\begin{equation}\label{eq:assoc-generalized-derivation}
T(ab)=T(a)b+aT(b)-a\xi b,
\qquad a,b\in A.
\end{equation}
Here and below, the module actions on $X^{**}$ are understood
to be the canonical extensions of the original module actions.

Put $z=T(1)$. Taking $a=b=1$ in
\eqref{eq:assoc-generalized-derivation}, and using the
unitality of the extended module actions, we obtain
\[
z=2z-\xi,
\]
and hence
\begin{equation}\label{eq:xi-z}
\xi=z.
\end{equation}
In particular, the correction element in
\eqref{eq:assoc-generalized-derivation} actually belongs to
the original module $X$.

On the other hand, since $\Psi_T(a,a)=0$ for every $a\in A$,
the definition of $K_T$ yields
\[
T(a^2)=2K_T(a)\circ a.
\]
Recalling that
\[
K_T(a)=T(a)-\frac14(az+za),
\]
we obtain
\begin{equation}\label{eq:jordan-square}
T(a^2)
=
aT(a)+T(a)a
-\frac14\bigl(a^2z+2aza+za^2\bigr).
\end{equation}
Meanwhile, taking $b=a$ in
\eqref{eq:assoc-generalized-derivation} and using
\eqref{eq:xi-z}, we get
\begin{equation}\label{eq:assoc-square}
T(a^2)=aT(a)+T(a)a-aza.
\end{equation}
Comparing \eqref{eq:jordan-square} and
\eqref{eq:assoc-square}, we conclude that
\begin{equation}\label{eq:commutator-square}
a^2z-2aza+za^2=0,
\qquad a\in A.
\end{equation}

We next prove that $z$ is central in $X$.
Fix a selfadjoint element $a\in A$, and consider the
bounded operators $L_a,R_a\in\mathcal B(X)$ defined by
\[
L_a(x)=ax,\qquad R_a(x)=xa.
\]
Since the left and right module actions commute, we have
\[
L_aR_a=R_aL_a.
\]
Set
\[
B=L_a-R_a.
\]
Equation \eqref{eq:commutator-square} gives
\[
B^2z=0.
\]
Therefore, for every $t\in\mathbb R$, the norm-convergent
exponential series implies
\[
e^{itB}z=z+itBz.
\]
Since $L_a$ and $R_a$ commute, we also have
\[
e^{itB}z
=
e^{itL_a}e^{-itR_a}z
=
e^{ita}ze^{-ita}.
\]
The elements $e^{ita}$ and $e^{-ita}$ are unitaries in $A$.
By contractivity of the module actions,
\[
\|e^{itB}z\|
=
\|e^{ita}ze^{-ita}\|
\leq \|z\|.
\]
Consequently,
\[
|t|\,\|Bz\|
=
\|itBz\|
\leq
\|z+itBz\|+\|z\|
\leq 2\|z\|.
\]
Letting $|t|\to\infty$, we obtain $Bz=0$, that is,
\[
az=za
\]
for every selfadjoint $a\in A$. Since every element of $A$
is a complex linear combination of selfadjoint elements,
it follows that
\begin{equation}\label{eq:z-central}
z\in Z_A(X).
\end{equation}

Now define
\[
d=T-L_z,\qquad d(a)=T(a)-az.
\]
Using \eqref{eq:assoc-generalized-derivation} and
\eqref{eq:xi-z}, we obtain, for all $a,b\in A$,
\[
\begin{aligned}
d(ab)
&=T(ab)-abz\\
&=T(a)b+aT(b)-azb-abz\\
&=(T(a)-az)b+a(T(b)-bz)\\
&=d(a)b+ad(b).
\end{aligned}
\]
Thus $d$ is a bounded derivation. Together with
\eqref{eq:z-central}, this proves that
\[
T=d+L_z,
\qquad d\in\operatorname{Der}(A,X),
\qquad z\in Z_A(X).
\]

It remains to prove the converse implication. Suppose that
\[
T=d+L_z,
\]
where $d$ is a bounded derivation and
$z\in Z_A(X)$.
Since every derivation on a unital algebra satisfies
$d(1)=0$, we have
\[
T(1)=z.
\]
Moreover, centrality of $z$ gives
\[
a\circ z=az=za.
\]
Consequently,
\[
K_T(a)
=
T(a)-\frac12(a\circ z)
=
d(a)+\frac12az.
\]
Using the Jordan derivation identity for $d$, we obtain
\[
\begin{aligned}
K_T(a)\circ b+a\circ K_T(b)
={}&d(a)\circ b+a\circ d(b)\\
&+\frac12(az)\circ b
+\frac12a\circ(bz)\\
={}&d(a\circ b)+(a\circ b)z\\
={}&T(a\circ b).
\end{aligned}
\]
In the second equality we have used the centrality of $z$,
which implies
\[
(az)\circ b=a\circ(bz)=(a\circ b)z.
\]
It follows that
\[
\Psi_T(a,b)=0,
\qquad a,b\in A.
\]

We have therefore established the equivalence between
$\Psi_T=0$ and the existence of a bounded derivation
$d$ and a central element $z\in X$ such that
$T=d+L_z$.

Finally, the representation is unique: evaluating
$T=d+L_z$ at the identity gives $z=T(1)$, and hence
$d=T-L_{T(1)}$.

\end{proof}

The zero-defect classification and a norm-stability estimate are distinct assertions. The following theorem identifies exactly what must be added to the former to obtain the latter.

\begin{theorem}\label{thm:closed-range}
For a unital contractive Banach $\A$-bimodule $X$, the space $\J(\A,X)$ is closed in $\B(\A,X)$. The following conditions are equivalent:
\begin{enumerate}
\item There is $C<\infty$ such that
\begin{equation}\label{eq:scalar-stability}
 \dist(T,\J(\A,X))\leq Cq(T)
 \quad(T\in\B(\A,X)).
\end{equation}
\item The operator $\mathscr R_1$ in \eqref{eq:scalar-R} has closed range.
\end{enumerate}
If its range is not closed, there are bounded linear maps $T_k$ with
\[
 \dist(T_k,\J(\A,X))=1,\qquad
 1\leq\norm{T_k}\leq2,\qquad
 \varepsilon_1(T_k)\longrightarrow0,
 \qquad q(T_k)\longrightarrow0.
\]
\end{theorem}

\begin{proof}
Proposition~\ref{prop:exact} identifies $\J(\A,X)$ with $\ker\mathscr R_1$. The operator is bounded by Lemma~\ref{lem:witness}; consequently its kernel is closed. The quotient
\[
 E=\B(\A,X)/\J(\A,X)
\]
is therefore a Banach space. The induced operator
\[
 \widetilde{\mathscr R}_1:E\longrightarrow\Bil_s(\A\times\A,X),
 \qquad [T]\longmapsto\Psi_T
\]
is bounded and injective. Its domain norm is
\[
 \norm{[T]}_E=\dist(T,\J(\A,X)).
\]

Assume (ii). The induced operator is a bounded bijection from $E$ onto the Banach space $\ran\mathscr R_1$. The bounded inverse theorem gives $C_0<\infty$ with
\[
 \dist(T,\J(\A,X))\leq C_0\norm{\Psi_T}
 \leq4C_0q(T).
\]
This proves (i).

Conversely, assume (i). Since $q(T)\leq\norm{\Psi_T}$, the induced operator is bounded below:
\[
 \norm{[T]}_E\leq C\norm{\widetilde{\mathscr R}_1[T]}.
\]
If a sequence $\Psi_{T_k}$ converges in the codomain, this inequality applied to differences shows that $[T_k]$ is Cauchy in $E$. It converges to some $[T]$, and continuity gives $\Psi_{T_k}\to\Psi_T$. Every norm limit of range elements is thus a range element, proving (ii).

Finally, if the range is not closed, the induced injective operator cannot be bounded below, by the preceding argument. For each $k$ choose $e_k\in E$ with $\norm{e_k}=1$ and $\norm{\widetilde{\mathscr R}_1e_k}<1/k$. By the definition of a quotient norm, $e_k$ has a representative $T_k$ of norm at most $2$. Its distance to the kernel is exactly $1$, so also $\norm{T_k}\geq1$. Changing a representative does not change its image under $\mathscr R_1$. Thus $\varepsilon_1(T_k)<1/k$ and $q(T_k)\leq\varepsilon_1(T_k)\to0$, as required.
\end{proof}

\begin{remark}
The last alternative is conditional on nonclosedness for a specified pair $(\A,X)$. It is a functional-analytic construction, not an assertion that a particular $C^*$-algebra and module exhibit failure of scalar stability. It also shows why exact rigidity alone is insufficient to infer a uniform quantitative theorem.
\end{remark}

\section{Matricial stability and approximate-diagonal correction}

From now on $X$ is a unital operator $\A$-bimodule in the concrete sense of Section~2. If $z=T(1)$, put $z_n=I_n\otimes z$. The canonical witness for $T_n$ is precisely $(K_T)_n$, since the $(i,j)$ entry of $A\circ z_n$ is $a_{ij}\circ z$. Define
\begin{align}
 \Psi^{(n)}_T(A,B)
 &=T_n(A\circ B)-(K_T)_n(A)\circ B-A\circ(K_T)_n(B),\label{eq:matrix-psi}\\
 \varepsilon_n(T)&=\norm{\Psi^{(n)}_T},
 &q_n(T)&=q(T_n),\label{eq:matrix-defects}\\
 \varepsilon_{\mat}(T)&=\sup_{n\geq1}\varepsilon_n(T),
 &q_{\mat}(T)&=\sup_{n\geq1}q_n(T).
\end{align}
The definition of $q(T_n)$ uses the algebra $M_n(\A)$ and module $M_n(X)$. Its witnesses are arbitrary maps at that level; compatibility between witnesses at different levels is unnecessary. If $T$ is completely bounded, all the displayed suprema are finite and
\begin{equation}\label{eq:matrix-witness}
 q_n(T)\leq\varepsilon_n(T)\leq4q_n(T),\qquad
 q_{\mat}(T)\leq\varepsilon_{\mat}(T)\leq4q_{\mat}(T)
 \leq4\norm T_{\cb}.
\end{equation}
These statements follow from Lemma~\ref{lem:witness} at each level. The map in \eqref{eq:matrix-psi} is formed using matrix multiplication; it is not obtained by merely applying a scalar symmetric bilinear defect entrywise.

For a bounded $D:\A\to X$ define the ordered Leibniz residual
\begin{equation}\label{eq:leibniz}
 R_D(a,b)=D(ab)-D(a)b-aD(b),\qquad r(D)=\norm{R_D}.
\end{equation}
Its multiplicative matrix amplification is
\[
 R_D^{[n]}(A,B)=D_n(AB)-D_n(A)B-AD_n(B).
\]
We use the notation
\[
 r_{\mat}(D)=\sup_n\norm{R_D^{[n]}},\qquad
 c(z)=\sup_{\norm a\leq1}\norm{za-az},\qquad
 c_{\mat}(z)=\sup_n\sup_{\norm A\leq1}\norm{z_nA-Az_n}.
\]
This specifies the amplification of the bilinear residual explicitly and avoids ambiguity between different notions of complete boundedness for bilinear maps.

\begin{lemma}\label{lem:block}
Let $T\in\B(\A,X)$, $z=T(1)$, and $D=T-L_z$. Then
\begin{equation}\label{eq:block-bounds}
 c(z)\leq4\varepsilon_2(T),\qquad
 r(D)\leq4\varepsilon_2(T).
\end{equation}
If $T(1)=0$, the second constant improves to $2$. For completely bounded $T$,
\begin{equation}\label{eq:block-cb}
 c_{\mat}(z)\leq4\varepsilon_{\mat}(T),\qquad
 r_{\mat}(D)\leq4\varepsilon_{\mat}(T),
\end{equation}
and again the second constant improves to $2$ when $T(1)=0$.
\end{lemma}

\begin{proof}
First define the nonsymmetric residual associated with the canonical witness:
\[
 E_T(a,b)=T(ab)-K_T(a)b-aK_T(b).
\]
For $a,b\in\A$, consider
\[
 U=\begin{pmatrix}0&a\\0&0\end{pmatrix},\qquad
 V=\begin{pmatrix}0&0\\b&0\end{pmatrix}.
\]
Their norms are $\norm a$ and $\norm b$. Their Jordan product is $\frac12\diag(ab,ba)$, and direct matrix multiplication gives
\begin{equation}\label{eq:ordered-block}
 \Psi^{(2)}_T(U,V)=\frac12
 \begin{pmatrix}
 E_T(a,b)&0\\
 0&E_T(b,a)
 \end{pmatrix}.
\end{equation}
The diagonal entries have this order because the second entry in $(K_T)_2(U)\circ V$ is $bK_T(a)/2$, and the second entry in $U\circ(K_T)_2(V)$ is $K_T(b)a/2$. Contractive compression of the first diagonal block proves
\begin{equation}\label{eq:ordered-estimate}
 \norm{E_T}\leq2\varepsilon_2(T).
\end{equation}

To obtain the correction commutator without losing the factor from \eqref{eq:ordered-block}, take
\[
 P=\begin{pmatrix}a&0\\0&-a\end{pmatrix},\qquad
 W=\begin{pmatrix}0&1\\0&0\end{pmatrix}.
\]
Here $P\circ W=0$ and $(K_T)_2(P)\circ W=0$, by linearity of $K_T$. Since $K_T(1)=z/2$, we obtain
\begin{equation}\label{eq:commutator-block}
 \Psi^{(2)}_T(P,W)=
 \begin{pmatrix}0&(za-az)/4\\0&0\end{pmatrix}.
\end{equation}
Its norm is $\norm{za-az}/4$. Taking the supremum over the unit ball proves the first estimate in \eqref{eq:block-bounds}.

Expanding the definition of $K_T$ gives
\[
 E_T(a,b)=T(ab)-T(a)b-aT(b)
 +\tfrac14(zab+abz+2azb).
\]
Meanwhile $D=T-L_z$ satisfies
\[
 R_D(a,b)=T(ab)-T(a)b-aT(b)+azb.
\]
Subtracting these identities yields
\begin{equation}\label{eq:bridge}
 E_T(a,b)=R_D(a,b)
 +\tfrac14\big((za-az)b+a(bz-zb)\big).
\end{equation}
Thus $r(D)\leq\norm{E_T}+\frac12c(z)\leq4\varepsilon_2(T)$. If $z=0$, the correction term vanishes, and \eqref{eq:ordered-estimate} gives the bound $2\varepsilon_2(T)$.

For the matrix-level assertions, fix $n$ and replace $a,b$ in the two block substitutions by $A,B\in M_n(\A)$, replace the identity by $I_n$, and replace $z$ by $z_n$. These substitutions belong to $M_2(M_n(\A))=M_{2n}(\A)$ and have the same norms as $A,B$. The preceding calculations then prove
\[
 \norm{E_T^{[n]}}\leq2\varepsilon_{2n}(T),\qquad
 \sup_{\norm A\leq1}\norm{z_nA-Az_n}
 \leq4\varepsilon_{2n}(T),
\]
where $E_T^{[n]}(A,B)=T_n(AB)-(K_T)_n(A)B-A(K_T)_n(B)$. Identity~\eqref{eq:bridge} holds with these matrix products as well. By \eqref{eq:matrix-actions}, it follows that
\[
 \norm{R_D^{[n]}}\leq4\varepsilon_{2n}(T),
\]
or at most $2\varepsilon_{2n}(T)$ if $z=0$. Taking suprema over $n$ proves \eqref{eq:block-cb}.
\end{proof}

\begin{corollary}\label{cor:block-exact}
For a unital operator bimodule, $\varepsilon_2(T)=0$ if and only if $T\in\J(\A,X)$. This conclusion follows from the block calculations without using the exact Jordan-to-associative theorem.
\end{corollary}

\begin{proof}
If $\varepsilon_2(T)=0$, Lemma~\ref{lem:block} gives $c(T(1))=0$ and $R_{T-L_{T(1)}}=0$. Therefore $T(1)$ is central and $T-L_{T(1)}$ is a bounded derivation. Conversely, if $T=d+L_z$ with $z$ central, then $d_2$ is a derivation for matrix multiplication and $z_2$ commutes with every matrix over $\A$. The converse calculation in Proposition~\ref{prop:exact}, applied at level two, gives $\varepsilon_2(T)=0$.
\end{proof}

The next lemma is the norm-controlled version of the usual degree-one approximate-diagonal argument. We prove both the derivation correction and the central correction, including the matrix estimates needed later.

\begin{lemma}\label{lem:correction}
Suppose $\A$ has a bounded approximate diagonal with bound $M$ and $X$ is a weak$^*$-closed unital operator $\A$-bimodule. Given $D\in\B(\A,X)$ and $z\in X$, there exist $\eta\in X$ and $z_0\in Z_{\A}(X)$ such that
\begin{equation}\label{eq:correction-ordinary}
 \norm\eta\leq M\norm D,\qquad
 \norm{D-d_\eta}\leq Mr(D),\qquad
 \norm{z-z_0}\leq Mc(z).
\end{equation}
If $D$ is completely bounded, the same elements can be chosen to satisfy
\begin{equation}\label{eq:correction-cb}
 \norm{D-d_\eta}_{\cb}\leq Mr_{\mat}(D).
\end{equation}
If an exact diagonal $m$ of norm at most $M$ exists, all conclusions hold for a norm-closed operator bimodule, without weak$^*$ closedness.
\end{lemma}

\begin{proof}
Define bounded linear maps on the projective tensor product by
\[
 F_D(u\otimes v)=uD(v),\qquad P_z(u\otimes v)=uzv.
\]
The defining norm of that tensor product and contractivity imply
\[
 \norm{F_D}\leq\norm D,\qquad \norm{P_z}\leq\norm z.
\]
Let
\[
 \eta_\lambda=F_D(m_\lambda),\qquad
 z_\lambda=P_z(m_\lambda).
\]
These nets lie in the balls of radii $M\norm D$ and $M\norm z$. The product of these two weak$^*$-compact balls is compact. Passing to a common subnet gives weak$^*$ limits $\eta,z_0\in X$, with $\norm\eta\leq M\norm D$. All subsequent limits are along this subnet.

For $a\in\A$, define a bounded tensor map
\[
 Q_{D,a}(u\otimes v)=uR_D(v,a).
\]
It satisfies $\norm{Q_{D,a}}\leq r(D)\norm a$. We claim that
\begin{equation}\label{eq:homotopy-identity}
 a\eta_\lambda-\eta_\lambda a
 =F_D(a\cdot m_\lambda-m_\lambda\cdot a)
 +\pi(m_\lambda)D(a)+Q_{D,a}(m_\lambda).
\end{equation}
For a simple tensor, this is exactly the identity
\[
 auD(v)-uD(v)a
 =auD(v)-uD(va)+uvD(a)+uR_D(v,a).
\]
It extends first by linearity and then by continuity to the completed tensor product, proving \eqref{eq:homotopy-identity}.

The first term on the right side of \eqref{eq:homotopy-identity} converges to zero in norm by approximate centrality. The second converges in norm to $D(a)$ because $\pi(m_\lambda)\to1$. The left side converges weak$^*$ to $d_\eta(a)$ by continuity of the module actions. It follows that
\[
 Q_{D,a}(m_\lambda)\longrightarrow d_\eta(a)-D(a)
 \quad\hbox{weak$^*$}.
\]
Each term has norm at most $Mr(D)\norm a$. Norm balls in $X$ are weak$^*$ closed, so
\[
 \norm{d_\eta(a)-D(a)}\leq Mr(D)\norm a.
\]
Taking the supremum proves the derivation estimate in \eqref{eq:correction-ordinary}. The map $d_\eta$ is an exact derivation because
\[
 ab\eta-\eta ab=(a\eta-\eta a)b+a(b\eta-\eta b).
\]

For centrality, the tensor identities give
\[
 az_\lambda-z_\lambda a
 =P_z(a\cdot m_\lambda-m_\lambda\cdot a)\longrightarrow0
 \quad\hbox{in norm}.
\]
Taking weak$^*$ limits proves $az_0=z_0a$ for all $a$, so $z_0\in Z_{\A}(X)$.
Define another bounded tensor map by
\[
 C_z(u\otimes v)=u(zv-vz).
\]
It has norm at most $c(z)$, and
\begin{equation}\label{eq:central-correction}
 z_\lambda-\pi(m_\lambda)z=C_z(m_\lambda).
\end{equation}
The left side converges weak$^*$ to $z_0-z$, and the right side has norm at most $Mc(z)$. Weak$^*$ closedness of the corresponding ball therefore proves the last bound in \eqref{eq:correction-ordinary}.

We now prove the completely bounded estimate using the same subnet and the same $\eta$. Fix $n$ and $A\in M_n(\A)$. Define the bounded tensor map with values in $M_n(X)$ by
\[
 Q_{D,A}^{(n)}(u\otimes v)
 =(I_n\otimes u)R_D^{[n]}(I_n\otimes v,A).
\]
Its norm is at most $r_{\mat}(D)\norm A$ by the matrix module inequalities. Its $(i,j)$ entry is $uR_D(v,a_{ij})$. The entrywise form of \eqref{eq:homotopy-identity} shows that
\[
 Q_{D,A}^{(n)}(m_\lambda)
 \longrightarrow (d_\eta)_n(A)-D_n(A)
 \quad\hbox{weak$^*$ in }M_n(X).
\]
For clarity, the entrywise centrality terms converge to zero in norm for each of the finitely many $a_{ij}$, and $\pi(m_\lambda)D(a_{ij})\to D(a_{ij})$ in norm. This justifies the displayed convergence without requiring uniform convergence over $n$ or over the unit ball.
Since every term has norm at most $Mr_{\mat}(D)\norm A$, weak$^*$ closedness of matrix norm balls gives
\[
 \norm{(d_\eta)_n(A)-D_n(A)}
 \leq Mr_{\mat}(D)\norm A.
\]
Taking the supremum first over $A$ and then over $n$ proves \eqref{eq:correction-cb}.

Finally, if $m$ is an exact diagonal, define $\eta=F_D(m)$ and $z_0=P_z(m)$ directly. The centrality terms in \eqref{eq:homotopy-identity} vanish, and $\pi(m)=1$. All the preceding estimates then follow as norm estimates of tensor maps, with no limiting argument. Both elements belong to the norm-closed module $X$, proving the final assertion.
\end{proof}

\begin{theorem}\label{thm:stability}
Suppose $\A$ has a bounded approximate diagonal with bound $M$, and $X$ is a weak$^*$-closed unital operator $\A$-bimodule.
\begin{enumerate}
\item For every $T\in\B(\A,X)$, there exist $\eta\in X$ and $z_0\in Z_{\A}(X)$ such that, writing $S=d_\eta+L_{z_0}$,
\begin{align}
 \norm{T-S}&\leq8M\varepsilon_2(T),\label{eq:main-ordinary}\\
 \norm{T(1)-z_0}&\leq4M\varepsilon_2(T),\label{eq:main-central}\\
 \norm\eta&\leq M\norm{T-L_{T(1)}}\leq2M\norm T.
\end{align}
In particular,
\begin{equation}\label{eq:q2-distance}
 \tfrac12q_2(T)\leq\dist(T,\J(\A,X))
 \leq32M q_2(T).
\end{equation}
\item For every $T\in\CB(\A,X)$, such $\eta,z_0$ can be chosen with
\begin{equation}\label{eq:main-cb}
 \norm{T-S}_{\cb}\leq8M\varepsilon_{\mat}(T).
\end{equation}
Consequently,
\begin{equation}\label{eq:qmat-distance}
 q_{\mat}(T)\leq\dist_{\cb}(T,\J_{\cb}(\A,X))
 \leq32M q_{\mat}(T).
\end{equation}
\item If $T(1)=0$, one may choose $z_0=0$ and improve \eqref{eq:main-ordinary} to $\norm{T-d_\eta}\leq2M\varepsilon_2(T)$. In the completely bounded case the estimate is $\norm{T-d_\eta}_{\cb}\leq2M\varepsilon_{\mat}(T)$.
\end{enumerate}
If $\A$ has an exact diagonal of norm at most $M$, weak$^*$ closedness of $X$ may be replaced by norm closedness throughout.
\end{theorem}

\begin{proof}
Set $z=T(1)$ and $D=T-L_z$. Lemma~\ref{lem:block} gives $r(D)\leq4\varepsilon_2(T)$ and $c(z)\leq4\varepsilon_2(T)$. Apply Lemma~\ref{lem:correction} simultaneously to $D,z$, obtaining $\eta,z_0$. Then
\[
 T-(d_\eta+L_{z_0})=(D-d_\eta)+L_{z-z_0}.
\]
Since $\norm{L_{z-z_0}}=\norm{z-z_0}$, the two correction estimates imply
\begin{align*}
 \norm{T-S}
 &\leq Mr(D)+Mc(z)\\
 &\leq8M\varepsilon_2(T).
\end{align*}
They also prove \eqref{eq:main-central} and the implementer estimate. The latter uses $\norm z\leq\norm T$ and $\norm{L_z}=\norm z$.

To obtain the upper inequality in \eqref{eq:q2-distance}, use $S\in\J(\A,X)$ and $\varepsilon_2(T)\leq4q_2(T)$. For the lower inequality, let $S_1\in\J(\A,X)$ be arbitrary. Its second amplification has exact canonical witnesses, by Corollary~\ref{cor:block-exact}. Retaining these witnesses for $T_2$ gives
\[
 q_2(T)\leq\norm{(T-S_1)_2}\leq2\norm{T-S_1}.
\]
The second estimate was proved in Section~2. Taking the infimum over $S_1$ proves the desired lower bound.

If $T$ is completely bounded, then so is $D$, because $L_z$ is completely bounded. Apply the matrix estimates in Lemma~\ref{lem:block} and the completely bounded part of Lemma~\ref{lem:correction}. The same decomposition now gives
\begin{align*}
 \norm{T-S}_{\cb}
 &\leq Mr_{\mat}(D)+Mc(z)\\
 &\leq Mr_{\mat}(D)+Mc_{\mat}(z)\\
 &\leq8M\varepsilon_{\mat}(T).
\end{align*}
The map $S$ is completely bounded, since it is the sum of an inner derivation and a module multiplier. This proves \eqref{eq:main-cb} and the upper inequality in \eqref{eq:qmat-distance}. For every $S_1\in\J_{\cb}(\A,X)$ and every $n$, its amplification is an exact third-type Jordan map, so
\[
 q_n(T)\leq\norm{(T-S_1)_n}\leq\norm{T-S_1}_{\cb}.
\]
Taking the supremum over $n$ and then the infimum over $S_1$ proves the lower inequality.

If $T(1)=0$, then $z=0$, $D=T$, and the central correction can be chosen as $z_0=0$. Lemma~\ref{lem:block} gives the sharper ordered-residual constants $2$, and Lemma~\ref{lem:correction} multiplies them by $M$. This proves (iii). The exact-diagonal assertion follows from the last part of Lemma~\ref{lem:correction} in each of the preceding steps.
\end{proof}

\begin{corollary}\label{cor:nuclear}
Let $\A$ be a unital nuclear $C^*$-algebra and $X$ a weak$^*$-closed unital operator $\A$-bimodule. For any bound $M$ on a bounded approximate diagonal of $\A$, all conclusions of Theorem~\ref{thm:stability} hold.
\end{corollary}

\begin{proof}
Haagerup's amenability theorem \cite{Haagerup}, together with the approximate-diagonal characterization \cite{JohnsonCohomology,JohnsonDiagonals}, supplies a bounded approximate diagonal. Its bound is finite. Substituting that bound into Theorem~\ref{thm:stability} proves the result.
\end{proof}

\begin{proposition}\label{prop:matrix-diagonal}
For $\A=M_d(\C)$, the tensor
\begin{equation}\label{eq:matrix-diagonal}
 m_d=\frac1d\sum_{i,j=1}^d e_{ij}\otimes e_{ji}
\end{equation}
is an exact diagonal and $\norm{m_d}_\pi=1$.
\end{proposition}

\begin{proof}
For a matrix unit $e_{kl}$, direct multiplication gives
\[
 e_{kl}\cdot m_d=\frac1d\sum_j e_{kj}\otimes e_{jl}
 =m_d\cdot e_{kl}.
\]
Linearity proves centrality for every matrix. Also
\[
 \pi(m_d)=\frac1d\sum_{i,j}e_{ij}e_{ji}
 =\frac1d\sum_{i,j}e_{ii}=1.
\]
These calculations establish the diagonal identities, but the displayed finite sum alone gives a poor projective norm bound. To obtain the optimal norm for this tensor, let $\mu$ be normalized Haar measure on the compact unitary group $U(d)$ and set
\[
 m=\int_{U(d)}u\otimes u^*\,d\mu(u).
\]
The integral is a Bochner integral in a finite-dimensional tensor space. Since $\norm u=\norm{u^*}=1$, $\norm m_\pi\leq1$.

We verify that this integral equals \eqref{eq:matrix-diagonal}. Invariance under left and right diagonal unitaries shows that
\[
 \int u_{ij}\overline{u_{lk}}\,d\mu(u)=0
 \quad\hbox{unless }i=l\hbox{ and }j=k.
\]
For fixed $j$, left permutation invariance makes $\int|u_{ij}|^2\,d\mu(u)$ independent of $i$. The column identity $\sum_i|u_{ij}|^2=1$ therefore makes this integral equal to $1/d$. Expanding $u\otimes u^*$ entry by entry gives $m=m_d$.
Thus $\norm{m_d}_\pi\leq1$. Conversely, the multiplication map is contractive, so
\[
 1=\norm{\pi(m_d)}\leq\norm{m_d}_\pi.
\]
The two inequalities prove $\norm{m_d}_\pi=1$.
\end{proof}

\begin{corollary}\label{cor:uniform}
Let $d\geq1$ and let $X$ be any norm-closed unital operator $M_d(\C)$-bimodule. For every bounded linear $T:M_d(\C)\to X$ there are $\eta\in X$ and $z_0\in Z_{M_d(\C)}(X)$ with
\[
 \norm{T-(d_\eta+L_{z_0})}\leq8\varepsilon_2(T)
 \leq32q_2(T).
\]
For completely bounded $T$,
\[
 \norm{T-(d_\eta+L_{z_0})}_{\cb}
 \leq8\varepsilon_{\mat}(T)\leq32q_{\mat}(T).
\]
If $T(1)=0$, the constants $8$ improve to $2$. None of these constants depends on $d$, $H$, $K$, or $X$.
\end{corollary}

\begin{proof}
Proposition~\ref{prop:matrix-diagonal} supplies an exact diagonal with bound $M=1$. The exact-diagonal version of Theorem~\ref{thm:stability} applies to the norm-closed module and gives every asserted estimate.
\end{proof}

\begin{remark}\label{rem:explicit}
The proof is constructive for $M_d(\C)$. With $z=T(1)$ and $D=T-L_z$, one can use
\[
 \eta=\frac1d\sum_{i,j}e_{ij}D(e_{ji}),\qquad
 z_0=\frac1d\sum_{i,j}e_{ij}ze_{ji}.
\]
The module products in these formulas use the specified representations. Lemma~\ref{lem:correction} proves the norm estimates for these elements by using the norm-one tensor representation, even though estimating the finite sums term by term would introduce a factor depending on $d$.
\end{remark}

\begin{corollary}\label{cor:matrix-closed}
Under the hypotheses of Theorem~\ref{thm:stability}, the operator
\[
 \mathscr R_2:\B(\A,X)\to
 \Bil_s(M_2(\A)\times M_2(\A),M_2(X)),
 \qquad T\mapsto\Psi_T^{(2)},
\]
has closed range. The operator
\[
 \mathscr R_{\mat}:\CB(\A,X)\to
 \ell^\infty\!\left(\Bil_s(M_n(\A)\times M_n(\A),M_n(X))\right),
 \qquad T\mapsto(\Psi_T^{(n)})_{n\geq1},
\]
also has closed range. The same assertions hold under the exact-diagonal hypothesis for norm-closed operator bimodules.
\end{corollary}

\begin{proof}
The first operator is bounded because $\varepsilon_2(T)\leq4\norm{T_2}\leq8\norm T$. Corollary~\ref{cor:block-exact} identifies its kernel with $\J(\A,X)$. Theorem~\ref{thm:stability} bounds the quotient norm by $8M\norm{\mathscr R_2(T)}$. The Cauchy-sequence argument in Theorem~\ref{thm:closed-range} therefore proves its range is closed.

The second operator is bounded because $\varepsilon_{\mat}(T)\leq4\norm T_{\cb}$. Its kernel is $\J_{\cb}(\A,X)$: vanishing at level two gives the required derivation and central element, and conversely their amplifications have zero defect at every level. The completely bounded estimate in Theorem~\ref{thm:stability} bounds the corresponding quotient norm by $8M$ times the range norm. The same complete-quotient argument proves closedness. The exact-diagonal case uses the identical estimates.
\end{proof}

\begin{example}\label{ex:noncentral}
Let $\A=X=M_2(\C)$ and $h=\diag(1,-1)$. Define $T(a)=ah$. It satisfies the exact associative generalized-derivation identity
\[
 T(ab)=T(a)b+aT(b)-ahb.
\]
Nevertheless,
\begin{equation}\label{eq:example-distance}
 \dist(T,\J(\A,\A))=
 \dist_{\cb}(T,\J_{\cb}(\A,\A))=1,
 \qquad \varepsilon_2(T)\geq\tfrac12.
\end{equation}
\end{example}

\begin{proof}
The associative identity follows by substituting $T(a)=ah$ and cancelling $ahb$. Any $S\in\J(\A,\A)$ satisfies $S(1)=\lambda1$ for some $\lambda\in\C$. Hence
\[
 \norm{T-S}\geq\norm{h-\lambda1}
 =\max\{|1-\lambda|,|-1-\lambda|\}\geq1,
\]
where the last inequality follows from $2\leq|1-\lambda|+|-1-\lambda|$. The same lower bound holds in completely bounded norm.

For the matching upper bound choose $\eta=h/2$ and $S=d_\eta$. Then
\[
 (T-S)(a)=\tfrac12(ah+ha).
\]
Both left and right multiplication by $h$ have completely bounded norm one, so $\norm{T-S}_{\cb}\leq1$. This proves the two equalities in \eqref{eq:example-distance}. Finally, $he_{12}-e_{12}h=2e_{12}$, and the commutator block \eqref{eq:commutator-block} gives $\varepsilon_2(T)\geq2/4=1/2$. Thus the central correction detected by the block test is a substantive part of the third-type Jordan problem.
\end{proof}

\section{Conclusion}

The third-type generalized Jordan identity admits a canonical normalization that is quantitatively comparable to arbitrary auxiliary witnesses. In a unital contractive Banach bimodule, its exact solutions are derivations plus central multipliers. At the scalar level, the resulting stability problem is equivalent to closedness of the range of a bounded linear defect operator. At matrix levels, two block substitutions recover the ordered Leibniz residual and the correction commutator. Combining these estimates with an approximate diagonal gives the constants $8M$ in operator and completely bounded norms, and $2M$ for maps vanishing at the identity. Full matrix algebras admit the dimension-independent constants $8$ and $2$.


\begin{thebibliography}{99}

\bibitem{ZeroProduct}
J.~Alaminos, J.~Extremera, M.~L.~C.~Godoy, and A.~R.~Villena,
\emph{Strongly zero product determined Banach algebras},
Linear Algebra Appl. \textbf{630} (2021), 326--354.
\href{https://doi.org/10.1016/j.laa.2021.09.002}{doi:10.1016/j.laa.2021.09.002}.

\bibitem{NearlyJordan}
M.~Eshaghi Gordji and N.~Ghobadipour,
\emph{Nearly generalized Jordan derivations},
Math. Slovaca \textbf{61} (2011), no.~1, 55--62.
\href{https://doi.org/10.2478/s12175-010-0059-x}{doi:10.2478/s12175-010-0059-x}.

\bibitem{Haagerup}
U.~Haagerup,
\emph{All nuclear $C^*$-algebras are amenable},
Invent. Math. \textbf{74} (1983), no.~2, 305--319.
\href{https://doi.org/10.1007/BF01394319}{doi:10.1007/BF01394319}.

\bibitem{HPS}
A.~Hosseini, A.~M.~Peralta, and S.~Su,
\emph{On the equivalence of all notions of generalized derivations whose domain is a $C^*$-algebra},
Linear Multilinear Algebra \textbf{73} (2025), no.~13, 3105--3124.
\href{https://doi.org/10.1080/03081087.2025.2499722}{doi:10.1080/03081087.2025.2499722}.
Preprint version consulted:
\href{https://arxiv.org/abs/2403.18773v2}{arXiv:2403.18773v2} (2024).
Numbered statements cited in the text refer to this preprint version.

\bibitem{JohnsonCohomology}
B.~E.~Johnson,
\emph{Cohomology in Banach algebras},
Mem. Amer. Math. Soc., no.~127, American Mathematical Society,
Providence, RI, 1972.

\bibitem{JohnsonDiagonals}
B.~E.~Johnson,
\emph{Approximate diagonals and cohomology of certain annihilator Banach algebras},
Amer. J. Math. \textbf{94} (1972), no.~3, 685--698.

\bibitem{JohnsonJordan}
B.~E.~Johnson,
\emph{Symmetric amenability and the nonexistence of Lie and Jordan derivations},
Math. Proc. Cambridge Philos. Soc. \textbf{120} (1996), no.~3, 455--473.
\href{https://doi.org/10.1017/S0305004100075010}{doi:10.1017/S0305004100075010}.

\bibitem{Ozawa}
N.~Ozawa,
\emph{Dixmier approximation and symmetric amenability for $C^*$-algebras},
J. Math. Sci. Univ. Tokyo \textbf{20} (2013), 349--374.
\href{https://arxiv.org/abs/1304.3523}{arXiv:1304.3523}.

\bibitem{Park}
C.~Park, S.~Ghaffary Ghaleh, and K.~Ghasemi,
\emph{Approximate $n$-Jordan $*$-derivations on $C^*$-algebras and $JC^*$-algebras},
J. Inequal. Appl. \textbf{2012} (2012), article~273.
\href{https://doi.org/10.1186/1029-242X-2012-273}{doi:10.1186/1029-242X-2012-273}.

\bibitem{PeraltaRusso}
A.~M.~Peralta and B.~Russo,
\emph{Automatic continuity of derivations on $C^*$-algebras and $JB^*$-triples},
J. Algebra \textbf{399} (2014), 960--977.
\href{https://doi.org/10.1016/j.jalgebra.2013.10.017}{doi:10.1016/j.jalgebra.2013.10.017}.

\bibitem{Ringrose}
J.~R.~Ringrose,
\emph{Automatic continuity of derivations of operator algebras},
J. London Math. Soc. (2) \textbf{5} (1972), no.~3, 432--438.
\href{https://doi.org/10.1112/jlms/s2-5.3.432}{doi:10.1112/jlms/s2-5.3.432}.

\bibitem{Sakai}
S.~Sakai,
\emph{On a conjecture of Kaplansky},
T\^ohoku Math. J. (2) \textbf{12} (1960), 31--33.

\bibitem{SameiFarsani}
E.~Samei and J.~Soltani Farsani,
\emph{Hyperreflexivity constants of the bounded $n$-cocycle spaces of group algebras and $C^*$-algebras},
J. Aust. Math. Soc. \textbf{109} (2020), no.~1, 112--130.
\href{https://doi.org/10.1017/S1446788719000089}{doi:10.1017/S1446788719000089}.

\end{thebibliography}
\end{document}